# A Study of the New Zealand Mathematics Curriculum

Neil Morrow (1), Elizabeth Rata (1) and Tanya Evans (2)

((1) School of Critical Studies in Education, University of Auckland, New Zealand,
(2) Department of Mathematics, University of Auckland, New Zealand)

Given the profound and uncritiqued changes that have been implemented in Aotearoa New Zealand education since the 1990s, this paper provides a critical commentary on the characterising features of the New Zealand mathematics' curriculum in the context of the first stage of a study. The emphasis is on the importance of research design that begins with an explicit, evidence-based hypothesis. To that end, we describe evidence that informs and identifies the study's hypothesised problem and causes. The study itself will show whether or not the hypothesis is justified; that is, is the absence of standardised prescribed content in New Zealand mathematics' curriculum the reason for the country's declining mathematics rankings? The study aims to increase understanding in the field of mathematics education by exploring the effects on New Zealand year 7 public school teachers' mathematics curriculum selection and design practices, teaching practices, and subsequently student achievement.

## Introduction

This paper describes the design of a research project which will examine the relationship between the New Zealand Mathematics curriculum and student achievement. Our purpose is to show the importance of the research design itself, particularly the first stage. It is in this stage that we hypothesise the existence of a problem with the New Zealand mathematics curriculum and decide on the methods of data collection and analysis for the empirical study. We also identify from the disciplinary literature the theoretical tools that contribute to the hypothesis of the problem's cause.

The aim of the research is to identify the features, causes and effects of the problem, so the first stage, the 'context of discovery', to use Reichenbach's term (Richardson, 2006) is to develop a justifiable hypothesis. Our focus on the hypothesis, which precedes the investigation itself (Reichenbach's 'context of justification') is to ensure that the enquiry itself is built on sufficient probability. According to an "objective Bayesian approach", a valid hypothesis will provide a sufficient degree of probability to the hypothesis. It "would not include a scientific hypothesis randomly chosen or irrational but would "put store on high degree of support of hypotheses by evidence" (Nola, 2001, p. 222). We would add that statements that assert the researcher's beliefs yet are without supporting evidence should not be included in a hypothesis. This can produce ideological manifestos rather than objective enquiries into social phenomena, a flaw prevalent in research of the advocacy type common to education studies (Bunge, 1998).

It is for this reason that we begin by identifying the problem in the New Zealand mathematics curriculum using evidence from international studies. This is followed by a section describing the proposed research methods to be used in an empirical study comparing New Zealand mathematics education to the Singapore system. The final section uses theories about knowledge, learning, and teaching to hypothesise the causes of the problem.

## The Research Problem

New Zealand's ongoing and long-term declining achievement in mathematics indicates a serious problem (Mullis et al., 2020; OECD, 2019). The country's underachievement in mathematics has been an issue for forty years, one first recognised in the Second International Mathematics Study (SIMS) results in 1981. The research showed that New Zealand's third form (year 9) students were performing at the lower quartile of participating countries (Hunter, 1998). The 1994/1995 the Third International Mathematics and Science Study found that the standard of mathematics learning in New Zealand was below the averages of fifty other countries in number (place value, fractions, and computation), measurement, and algebra at the age levels tested (Holton, 2005; Hughes & Peterson, 2003; Thomas & Tagg, 2005).

As a result of these findings, the government committed $75 million in an attempt to remedy the perceived inadequacies of teachers in mathematics, primarily the teachers' lack of mathematical knowledge, which was seen to be the cause of the problem. This led the Ministry of Education to introduce the Numeracy Development Project (NDP) from 2000 (Young-Loveridge, 2010). The NDP was part of a complete review of New Zealand's curriculum, which began in the 1990s and emerged in the completed form with the 2007 National Curriculum (Ministry of Education, 2007).

New Zealand's average mathematics achievement has continued to deteriorate since the new curriculum's introduction in 2007. Apart from Year 5's average achievement in the Trends in International Mathematics and Science Study (TIMMS), New Zealand's average (mean) achievement in mathematics has continued to decline in both the TIMMS and Programme for International Student Assessment (PISA) (Mullis et al., 2016; OECD, 2019). At both levels measured by TIMMS (Years 5 and 9), New Zealand's average achievement is still significantly lower than the centre point of the relevant TIMMS scale (Mullis et al., 2016). Given these large-scale assessments test different constructs of mathematics, TIMMS being curriculum oriented (pure mathematics), and PISA measuring the application of skills to real-life contexts (applied mathematics), New Zealand's average achievement in both pure and applied areas of mathematics is diminishing (Gronmo & Olson, 2006; Mullis et al., 2016; OECD, 2019). Since New Zealand's participation in PISA's first assessment of global mathematics achievement in 2003, not only has average performance declined, so too has achievement amongst the students at 'top and the bottom of the performance distribution' (OECD, 2019, p.321). Accordingly, only 12% of New Zealand's 15-year olds scored at the top two levels (Level 5 or higher) in mathematics compared to Singapore's 37% in the latest PISA cycle. Whereas at the bottom end of the performance distribution, 22% of New Zealand's cohort are 'low achievers in maths' contrasted to 2% of the students assessed in China (OECD, 2019). In other words, one in five of New Zealand's 15-year-olds are considered low achieving in mathematics.

Rather than address the need for teachers to improve their mathematical knowledge (and that of other subjects), which was the intention of changes to mathematics following the 1994/1995 Third International Mathematics and Science Study Report, the new outcomes-based curriculum was a radical break from the past with respect to knowledge itself. Instead of teachers being required to know their subject more thoroughly, by the time the finalised national curriculum was published (MoE, 2007), the problem of inadequate teacher knowledge was no longer a problem. This sleight of hand was achieved by removing the subject content itself. The standardised curriculum with its various prescribed syllabuses for academic subjects was replaced by a 'curriculum without content'.

The 2007 curriculum is a framework only. It contains a vision and values statement, key competencies, and general education principles. It is left to schools and teachers to decide what content to teach. Achievement Objectives listed at the back of the national curriculum document describe what students will be able to do – their skills and competencies – but the knowledge itself is not prescribed. The shift from content to competencies took from the early 1990s to the finalised document in 2007 to accomplish. Achievement Objectives had appeared in the mathematics curriculum published in 1992 (Learning Media, 1992), but this document also contained a considerable amount of specified content for all school levels – nearly 200 pages of detail for Number, Measurement (and calculus from level 7), Geometry, Algebra, Statistics, and Calculus. That specified content was not in the 2007 finalised curriculum.

The emptying out of mathematical knowledge (and of other subject knowledge) in the New Zealand Curriculum (MoE, 2007) is characterised by the shift from 'knowledge-that' (subject concepts and associated content) to 'know-how-to' (competencies and skills) (Rata, 2021; Ryle, 1949). This removal of knowledge itself was justified by another major change – to innatist or developmental learning (Winch, 1997), now seen in inquiry and personalised learning approaches from year one (Zame, 2019) to senior schooling. The 'learnification' (Biesta, 2012) of New Zealand education concealed the absence of the actual knowledge which had occurred in the shift to the 'know-how-to' of outcomes-based education. It also enabled a profound change to the teacher's role – from subject expert to facilitator of learning (Rata, 2017). The student, too, was re-imagined as 'the learner', one whose interests, background, and culture were now at the centre of schooling (for example, see the development of 'culturally responsive pedagogies' MoE, 2019).

The combination of outcomes-based education with the learning approach has cemented the fundamental change to New Zealand education precipitated by the removal of 'knowledge-that' (generalising concepts and materialised content) from the national curriculum in the post-1990s' decades. In the final section, we provide a more detailed discussion of the knowledge and learning theories that will be used to critique the 2007 curriculum. This includes looking more closely at the significance of the distinction between the two forms of propositional (i.e. academic) knowledge, cognitive theory's critique of innatist learning theory, and the role of the teacher as an autonomous curriculum maker.

## Research Methods

The proposed project is a comparative study of New Zealand's localised, generic, outcomes-based mathematics curriculum and Singapore's centrally prescribed national mathematics curriculum. The aim is to identify the differences between the two curriculum types and their effects on teaching methods and student achievement. Mixed methods will be employed to combine the complementary powers of qualitative and quantitative analyses.

## Data Sources

The study will collect comparative data from three sources: 1. the New Zealand and the Singapore mathematics curricula content and achievement statements), 2. The teaching methods, resources, and activities used to teach algebraic expressions and equations by mathematics teachers in each of the two countries respectively (matched sample), and 3. Data about student achievement related to that algebra topic.

### Content statements

A comparison will be made between the two countries' content statements. Levels 3-5 of the mathematics operational guidance section in the *NZC,* i.e. achievement objectives (MoE, 2007) will be compared to the Primary Six year level of Singapore's Mathematics Syllabus: Primary One to Six (Ministry of Education, 2012, 2020). Each curriculum document's remaining elements related to mathematics and/or generic features relating to every curriculum area will also provide data material. This includes the statement of purpose of mathematics in the middle of the *NZC* and statements at the 'front end' (Priestley & Sinnema, 2014), which describe the directions for learning for all learning areas (vision, values, key competencies, and principles). As the Singapore Mathematics Syllabus: Primary One to Six is exclusively a mathematics curriculum document, the disciplinary content, the syllabuses' specific pedagogical and philosophical statements, and those statements related to Singapore's mathematics curriculum will provide data for that country.

### Teacher instruction

The data collected about teaching practice will be from interviews, observations, and resource analysis of the teaching methods used by Year 7 New Zealand mathematics teachers and Primary Six Singapore mathematics teachers as they teach the same algebra topic: algebraic expressions and equations.

The interviews with the New Zealand and Singapore teachers, both pre-observation and post-observation, will cover a range of detail, including biographical background information such as their mathematical and teaching qualifications, teaching experience, beliefs about mathematics education, their own evaluation of their mathematical knowledge, teaching expertise, and decision-making ability. The interviews will also ask about the teachers' views concerning curriculum prescription and teacher autonomy. There will be questions specific to algebraic expressions and equations– about content selection, design and resources as well as teaching methods and learning activities. The teachers' planning documents (schemes of work, unit/topic plans, lesson plans, reference materials, assessment methods) will be used as baseline data when observing the actual teaching of the topic.

Those observations by the researcher of the teachers' actual practice will provide the third means of data collection. Field notes and video devices will record the teachers' instructional methods, the content used, the specific mathematical language, and the students' activities at the start, middle and end of the algebra topic over three weeks.

### Student achievement

The student achievement data will be collected using a topic assessment at the start and the end of the algebra topic unit, which will include examples of both knowledge types. The topic involves linear expressions or equations (the highest exponent of varibles is 1).

The assessment tool will involve distinct sections based on content that progresses from the identification and definition of concepts through to interpreting, representing and solving algebraic equations. Each content type involving procedures (operations) in the assessment will incorporate a hierarchal sequence from least complex to the most complex in relation to each specific content element.

## Data Analysis

Data analysis provides the meeting ground between the empirically obtained data and the use of sociological and epistemological concepts that will be used to explain what the data means. This explanatory work requires a theory of knowledge. Briefly, mathematical knowledge is propositional knowledge. It is theoretical, abstract, and conceptual. Although a mental construct, mathematical knowledge can be suggested by and applied to the material world. These two forms of mathematics – the theoretical and the applied – are referred to as the knowledge forms 'knowledge-that' (concepts and content) and 'know-how-to' (skills and competencies) (Rata, 2021; Ryle, 1949). We say more about knowledge theory and its usefulness as an analytical and explanatory tool in the following section.

### The analytical device

Our purpose here is to describe how the two knowledge forms will provide the main analytical categories in making sense of the data. All the data, including the broader elements of the curriculum documents, the material collected from the interviews, documents and observations, will be analysed in terms of the two main propositional knowledge forms – 'knowledge-that' and 'know-how-to'. The tool for these analyses will be the Curriculum Design Coherence (CDC) Model (Rata, 2021) because it 'makes explicit the interdependent relationship between the two knowledge forms, 'knowledge-that' (subject concepts and subject content) and 'know-how-to' (subject competencies, techniques, and skills) (Rata & McPhail, 2020, p.1). The CDC Model's differentiation between the two forms of propositional knowledge means that it will enable a thorough examination of the data to explore how those two forms are understood in the curriculum documents themselves and also by the teachers in their practice.

## Teachers data

The CDC Model will be applied to analysing three sets of data. The first set is the teachers' data – collected from the interviews, observations and teaching materials as they plan and teach the algebra topic. This will enable an analysis of the weighting given to one or other of the knowledge forms and the use of specific mathematical language. This will lead to explanations of the effects of this weighting on teacher knowledge selection and student mathematical understanding and performance.

The teachers' planning materials will be analysed to identify the selected content used to materialise the mathematical concepts of the algebra topic. The analysis will examine whether the selected content best 'demonstrates the concept's breadth and depth of abstraction' (Rata, 2016, p. 173; see also CDC Model Element 2, Criterion one in Rata, 2021, p. 18). The planning documents analysis will also enable the researcher to identify what informs the teachers' curriculum decision making, including such factors as resource selection and the way lessons are structured and delivered. This analysis will examine whether or not the teachers' curriculum designs provide a distinction and connection between subject concepts/content ('knowledge-that' knowledge) and competencies ('know-how-to knowledge).

The data collected from the observations of lessons will provide information about the language used in mathematics lessons and show the alignment between the intended and enacted curriculum – do the teachers actually do what they plan?

## Curriculum document analysis

The second use of the CDC Model will be in analysing the disciplinary content in the New Zealand and Singapore content and achievement statements. The analysis will look specifically at the extent to which the curriculum topics have a cumulative nature that promotes epistemic ascent (Winch, 2013) or conceptual progression (Rata, 2016). The hierarchal nature of mathematics means acquiring new concepts relies on students understanding the previous concepts in the progression-abstraction chain (Muller, 2006). This analysis will enable a comparison of the nature of content progression between the two distinct curriculum structures; the level-based structure utilised in the New Zealand mathematics curriculum and the year based spiral approach employed in Singapore's.

## Student achievement

The achievement data collected before the start and at the end of the algebra topic (algebraic expressions and equations) will be used in a repeated-measures statistical test to measure each student's progress while controlling for the following variables: content and knowledge types ('knowledge that' vs 'knowledge how') and the location (NZ vs Singapore). The use of the algebra topic in the classroom study is apposite because it is the first instance in both curriculums that students are introduced to algebra involving variables (Ministry of Education, 2012; MoE, 2007). Subsequently, this will affect a conceptual move from thinking about specific numbers and measures as found in arithmetic and measurement 'towards relations among sets of numbers and quantities, especially functional relations' (Carraher, Schliemann, Brizuela, & Earnest, 2006, p. 88).

# The New Zealand Curriculum

The hypothesis of the proposed study is that the ongoing decline of student mathematical achievement is the result of four main interdependent features which characterise the 2007 New Zealand curriculum. The features are: 1. Its highly generic non-prescriptive nature, 2. A commitment to teacher autonomy in curriculum knowledge selection, 3. A competency-based outcomes approach, 4. A commitment to localisation in curriculum selection. We discuss each one in turn to show how they contribute to and draw from the others to create a 'curriculum without content'. Of course, these characteristics have emerged from deeper socio-political forces and the ideological and intellectual ideas associated with those forces. Space prevents an in-depth account of the origins of each feature and the reasons for their convergence (see

Rata, 2012; 2021 for that discussion). However, it is essential to know about these underlying forces and influences in order to understand the extent to which mathematical knowledge, along with other academic subjects, has been affected in New Zealand education.

## A curriculum without content

The first feature is the curriculum's highly generic and non-prescriptive nature. The National Curriculum (MoE, 2007) 'sets the direction for teaching and learning … it is a framework rather than a detailed plan' and 'schools have considerable flexibility when determining the detail' (p. 37). That detail is the actual content – what is to be taught. Its absence from the national document means that there is no standardisation across the country. Individual schools and teachers select the mathematical content for their respective class and school. The mathematics taught in a school in one part of the country or in one part of a city may bear no relation to what is taught in another part.

In contrast, consistently high performing Singapore has a centrally planned curriculum that incorporates a single mathematics curriculum framework that shares a common emphasis throughout the levels (Ministry of Education, 2012; MoE, 2007; Soh, 2008). This framework unifies the direction of the mathematics curriculum for all levels from primary to pre-tertiary (Ministry of Education, 2012). There is a series of connected syllabuses, each with its specific aims designed to meet the 'different needs and abilities of students' (Ministry of Education, 2012, p. 9). Soh (2008) explains the rationale for the centrally planned nature of Singapore's curriculum:

> what a child needs to learn in mathematics in the formative years is common and requires careful thought and planning to make it accessible to every student. A centrally planned curriculum provides clear guidance in teaching and learning to teachers (p. 27).

The hierarchical sequencing enabled by its year-based curriculum design and the specific subject content is prescribed in the Content and Learning Experiences by Level element of each syllabus (Ministry of Education, 2012).

An example of the difference between the New Zealand and the Singapore curricula is found in the different approaches to whole numbers. The New Zealand curriculum has the generic Achievement Objectives at Level 1 'Use a range of counting, grouping, and equal-sharing strategies with whole numbers and fractions', 'Know the forward and backward counting sequences of whole numbers to 100', and 'Know groupings with five, within ten, and with ten' (MoE, 2007). In contrast, the Singapore Primary One (Year 1) Whole Numbers sub-strand of the Number and Algebra strand provides specific detail about subject content to be taught. The sub-strand comprises three sections: Numbers up to 100, Addition and Subtraction, Multiplication and Division. The content statements in each section range from five (Multiplication and Division) to eight (Numbers up to 100 and Addition and Subtraction). The first section, Numbers up to 100, has eight specific content statements:

**1. Numbers up to 100**
1.1 counting to tell the number of objects in a given set
1.2 number notation, representations, and place values (tens and ones)
1.3 reading and writing numbers in numeral and in words
1.4 comparing and ordering number of objects in two or more sets
1.5 comparing and ordering numbers
1.6 patterns in number sequences
1.7 ordinal numbers (first, second, up to tenth) and symbols (1st, 2nd, 3rd, etc)
1.8 number bonds for numbers up to 10
 (Ministry of Education, 2012, p. 34).

These content statements are supplemented by a fine-grained description of the learning experiences, which provide explicit curriculum information, guiding the teachers in the selection of knowledge:

**Students should have the opportunities to:**
(a) use number-bond posters and make number stories to build and consolidate number bonds for numbers up to 10.
(b) work in groups using concrete objects to
  • make a group of ten and count on from 10 to tell the number (less than 20).
  • make groups of ten and count tens and ones to tell the number (more than 20).
  • estimate the number of objects in a set before counting.
  • make sense of the size of 100.
(c) use concrete objects and the base-ten set to represent and compare numbers in terms of tens and ones, and use language such as 'more than', 'fewer than', 'the same as' and 'as many as' to describe the comparison.
(d) play games using dot cards, picture cards, numeral cards and number-word cards etc. for number recognition and comparison.
(e) describe a given number pattern using language such as '1 more/less' or '10 more/less before continuing the pattern or finding the missing number(s).
  (Ministry of Education, 2012, p. 34).

Furthermore, their implementation is supported by centrally authorised textbooks that "help teachers understand the emphases and scope of syllabuses" (Kaur et al., 2015, p. 313). At the primary level, instructional or pedagogical guides include thorough schemes of work that have lesson plans, common misconceptions, and show conceptual progression (Dindyal, 2006; Kaur et al., 2015).

## Teacher autonomy

The second feature which characterises the New Zealand curriculum and which is, we hypothesise, a contributing cause to the curriculum's failure is the commitment to 'teacher autonomy'. In the absence of prescribed mathematical content, New Zealand teachers have considerable autonomy in selecting what they will teach. Content is selected and designed from a wide range of resources. It may or may not include resources provided by a Ministry of Education website NZMaths (MoE, 2021). The material may come from an increasing number of commercial businesses, from social media websites such as Facebook, YouTube, Google, or from material the teacher has acquired before the new curriculum drafts from the 1990s.

The 2007 curriculum gave New Zealand teachers an autonomy over curriculum knowledge selection; an autonomy that had traditionally been exercised in the domain of pedagogy. Indeed, the country's progressive pedagogy with its creative teaching methods and teacher autonomy can be seen in the open-air movement of the 1920s and the New Education Fellowship of the 1930s (Couch, 2014). However, prior to 2007, teacher autonomy over pedagogical matters was counter-balanced by a national curriculum that specified content in various syllabus booklets for each subject.

The extension of teacher autonomy from teaching methods to include curriculum selection was the result of a profound change to how knowledge itself is understood. This led to the emptying out of content from the curriculum and to the conflation of curriculum and pedagogy – of the 'what' and the 'how' (Young, 2010). We discuss the knowledge issue below. A cause of direct relevance to the insistence on teacher autonomy was the post-1980s shift to teaching as a profession. Accordingly, teachers, especially primary school and early childhood teachers graduating from the new university-based teacher education faculties with bachelor degrees in teaching, were to have the same level of autonomy over their work which characterises other professions. Because the strengthening of the professional nature of teaching occurred in tension with an equally strong trend – that of managerialism (Benade, 2011), the result has tended to be an overly-zealous commitment to teacher professional autonomy, but now in its new form as autonomy over curriculum selection as well as the traditional autonomy with respect to pedagogy.

However, the arguments commonly used to justify teacher autonomy over curriculum selection are not based on rigorous evidence-based research. For example, in addressing critically important equity considerations in Aotearoa New Zealand, the work of Korthagen (2017) is often used to promote localised curriculum to enable culturally responsive and sustained pedagogies that take an asset-based approach (Celedón-Pattichis, Borden, Pape, Clements, Peters, Males, Chapman, & Leonard, 2018). The kernel of the argument is formed by extracting the following statement from Korthagen (2017): "attempts at influencing teacher behaviour have to be adjusted to individual teachers in their specific circumstances and settings, and that it is impossible to promote change through a pre-planned, fixed curriculum" (p. 391)." Under a closer inspection of the manuscript, the following realisation manifests acutely. At least three reasons make this quote unsuitable for inclusion as part of an argument to promote a culturally responsive and sustained curriculum.

First, the choice of the quote stands defenceless under the accusation of cherry-picking. The paper is about "a critical analysis of traditional and new approaches to professional development" (p. 387) and is not an empirical investigation of the impacts of different curricula. The quote is a side note made by the author, expressing his opinion that is not based on empirical findings and should be treated as a speculation. Second, the paper is not a scoping or a systematic review of the literature – it is a very narrow summary of research on the topic of teacher learning. The scientific rigour of this approach is questionable. The author states: "after discussing a framework on teacher learning, an approach will be presented representing a more radical attempt at integrating practice and theory, namely by giving the person of the teacher a more central place. In this approach, which I call *professional development 3.0*, the professional and the personal aspects of teaching are intertwined." (p. 389). It transpires that the paper is a collection of thoughts and opinions by the author presented as a convincing narrative to promote his new "more radical" approach to professional development. This is not a research study that involved data collection and rigorous analysis that withstands scientific scrutiny; it is not even a rigorously conducted review paper. Nevertheless, the author puts forward unsubstantiated claims that promise to solve all problems and alludes to the value of a non-planned curriculum. This is not an evidence-based recommendation.

Finally, even though the paper is not based on any empirical investigation, it contains a section titled "Empirical evidence" (p. 397), which is simply a summary of 'evidence' from different studies that are rather inappropriately used to support the main outcome of the study – a 'more radical' professional development framework. Most of the studies used as 'evidence' are small case studies - as small as studying 1 participant - in the Netherlands and the US. This is not research that can be taken seriously and utilised as a generalisable conclusion in the Aotearoa New Zealand context.

This type of theoretical deliberations that are not validated in empirical studies cannot produce generalisable findings and, thus, should never be used in policy and practice recommendations to avoid the inevitable failure that would result. It is likely that such recommendations will be ultimately blamed for policy failures in the court of history. Nevertheless, it seems that the proliferation of such research over the years has resulted in a self-referencing self-sustaining research paradigm proclaiming teacher autonomy over curriculum selection as an unquestionable 'right'.

However, teacher autonomy to select *what* is to be taught as well as *how* to teach it imposes impossible demands upon teachers. They are curriculum-makers, curriculum-designers, and curriculum-implementers. According to the findings from the Knowledge-Rich School Project (Rata, 2021) this overload leads to confusion and difficulty with both knowledge selection and design and a tendency to focus solely on the 'how' of teaching rather than the content itself.

## Curriculum Knowledge

The third major feature characterising the New Zealand curriculum and contributing to its limitations is the shift from prescribed academic knowledge of specified subjects such as mathematics to the outcomes approach and its accompanying focus on the 'learner' and inquiry learning. It is a shift that occurred to a greater or lesser extent in most developed countries from the 1990s. However, as with the 'learnification' movement (Biesta, 2012), competency-outcomes education was most enthusiastically taken up in New Zealand (Priestley & Sinnema, 2014). Spread by the OECD, though originally promoted by global corporations (Delors, 1996), Gilbert's book 'Catching the Knowledge Wave' publication was the seminal influence in this country.

The shift is fed by several contributing movements that were brought together in a major 2012 Ministry of Education commissioned report (Bolstad et. al., 2012) promoting what became known as '21st Century Education'. All aspects of education were affected including school architecture, along with fundamental changes to pedagogy and the curriculum. While space prevents an analysis of the reasons for the rapid and uncritical embrace of '21$^{st}$ century education' by the Ministry it is significant that the change to the role of education occurred at the height of corporate globalisation. Fundamental changes to production and distribution, and significantly for education, to the management of labour, promoted a restricted view of education's role in producing the flexible and highly skilled worker with dispositions, competencies, and the technological skills with which to locate 'just in time' knowledge (Gilbert, 2005). A documented account that tracks the 21$^{st}$ century competency movement from the OECD to, and then within, New Zealand is available in Lourie, 2020.

The shift from knowledge to competencies and outcomes occurred alongside and contributed to the 'learnification' movement (Biesta, 2012). The belief that 'learning' would arise from within the child through a curiosity stimulated by a facilitator using personalised inquiry pedagogies. Indeed, a study by Zame (2019) of inquiry learning by teachers of year one students shows the extent of the belief in knowledge being 'constructed' by the children as they are stimulated to inquiry. Zame's study also shows how the 'learning' that is claimed to occur does not in fact happen.

The extent to which the word 'learner' is used in official education documents demonstrates the degree to which this belief is normalised. However, innatist or development theories (Winch, 1998) are strongly criticised by cognitive scientists and evolutionary theorists (Geary, 2002; Geary & Berch, 2016). Cognitive and evolutionary theory distinguishes between academic knowledge (with its propositional character and disciplinary origin) and knowledge from experience (i.e. socio-cultural knowledge and proto-science). An account of the theories and their implications for education are available in Rata (2021).

The dominance of innatist learning theory (with its Rousseauian roots) is supported by two factors. The first is the absence of cognitive and evolutionary theories of learning mentioned above in Ministry policies concerned with learning. The second is the absence of a theory of knowledge with which to inform curriculum development. Indeed, recent drafts of proposed changes to science, history and mathematics appear informed by a mismatch of relativist beliefs about knowledge. There is no coherent theory of knowledge that recognises the distinction between knowledge from experience and the propositional knowledge, which forms the concepts and content of curriculum subjects (Cuthbert & Standish, 2021).

## Localisation

The fourth feature of the New Zealand curriculum follows on from the absence of a theory of knowledge. The localised curriculum (MoE, 2019) promotes knowledge, not as the objective scientifically verifiable truth, produced in the disciplines and altered for teaching at school as academic subjects, but as the experience-derived truth of local social groups. The conflation of two different types of knowledge is informed by the intellectual shift from universal knowledge to the knowledge relativism celebrated by postmodern writers and embraced by various identity politics movements, including those in New Zealand (Rata, 2017). In education, the country's wider commitment to biculturalism takes the form of an emphasis on cultural recognition in pedagogical relations, an approach strongly influenced by Bishop's seminal book Culture Counts (1999). The result is that socio-cultural knowledge and proto-science is added to the mix of the competencies and learning approaches.

Here are two illustrative examples. The first is from a study of Māori teachers' classroom practices (Lynch, 2017). The teachers had benefitted from an academic education themselves and intended this for their own children. However, in line with bicultural policy, they teach a socio-cultural curriculum to their Māori students. The social studies teacher has replaced history and geography with kapa haka to "provide students with an opportunity to learn … through a Māori lens" (p. 56). Another teacher rejected the idea of educational success, calling it "white success" and in opposition to succeeding "as Māori" (p. 60). The second example is from the media (Collins, 2020). According to the principal of a 21[st] century secondary school, the "dangers of prescribing a powerful knowledge curriculum" are because it "is about whose knowledge". A "Eurocentric" approach is "a colonial tool of putting old western knowledge ahead of indigenous communities".

The absence of a theory of knowledge that might then inform firstly what academic knowledge actually is, then how it should be selected, designed and taught contributes to one of the main forces behind the localisation of the curriculum. The examples above show that the shift to 'culture' affected not only pedagogy but had a major influence on what knowledge was to be selected by the 'autonomous' teacher. Although an emergent 'knowledge in education' movement (Barrett, Hoadley & Morgan, 2017) has contributed to a re-commitment to a knowledge curriculum in England from 2012, it appears to have gained little traction in New Zealand. This may be the result of two reasons. The first is the absence of a theory of knowledge from recent discussions about the curriculum. The second is that no distinction is made between types of knowledge, that is between socio-cultural knowledge and proto-science on the one hand and disciplinary (propositional) knowledge on the other (Bishop, 2019).

The combination of the four features we identify above suggests a context within which mathematics education is severely compromised. The extent to which each of these features contributes to the decline of student achievement will be investigated.

## Conclusion

This study will bring increased understanding to the field of mathematics education by exploring the effects of the generic nature of the national curriculum on New Zealand year 7 public school teachers' mathematics curriculum selection and design practices, teaching practices, and subsequently student achievement. In other words, we want to know 'what' knowledge New Zealand teachers are selecting, 'how' they are using that knowledge to design programmes in mathematics, and how this affects teaching practice and, in turn, student achievement. However, we wish to avoid the tendency in educational research to 'fit' a study to a set of ideological beliefs. Given the huge and uncritiqued changes that have been implemented in New Zealand education since the 1990s (and not only to the curriculum), our purpose in this paper is to show the importance of research design (the context of discovery) which begins with an explicit, evidence-based hypothesis. We have described evidence that identifies the study's hypothesised problem and causes. The study itself will show whether or not the hypothesis is justified; that is, is the absence of standardised prescribed content in New Zealand mathematics' curriculum the reason for the country's declining mathematics rankings?